\documentclass[a4paper]{article}

\usepackage{type1cm}   
\usepackage{dsfont}         
\usepackage{makeidx}         % allows index generation
\usepackage{graphicx}        % standard LaTeX graphics tool
\usepackage{multicol}        % used for the two-column index
\usepackage[bottom]{footmisc}% places footnotes at page bottom

\usepackage{newtxtext}       % 
\usepackage{newtxmath}       % selects Times Roman as basic font
\usepackage{pict2e}

\usepackage{algorithm}% http://ctan.org/pkg/algorithms
\usepackage{algpseudocode}
\algtext*{EndIf}
\algtext*{EndFor}
\algtext*{EndWhile}
\usepackage{amsmath}
\usepackage{subcaption}
\usepackage{mathrsfs}
\usepackage{tikz}
\usepackage{mathbbol}

\newtheorem{wdef}{Definition}

\newcommand{\wabs}[1]{{\left| {#1} \right|}}

\newcommand{\wcal}[1]{{\cal {#1}}}

\newcommand{\wds}[1]{{\mathds #1}}

\newcommand{\wdf}[2]{{#1}' \left(#2\right)}

\newcommand{\wdfc}[2]{{#1}'\!\left(#2\right)}

\newcommand{\wdsfc}[2]{{#1}''\!\left(#2\right)}

\newcommand{\wfc}[2]{{#1}\!\left(#2\right)}

\newcommand{\whs}[1]{\hspace{#1cm}}

\newcommand{\wit}[1]{\mathbb{#1}}

\newcommand{\wqes}{\hfill \ensuremath{\blacktriangle}}

\newcommand{\wlr}[1]{\left({#1}\right)}

\newcommand{\wo}[1]{\overline{#1}}

\newcommand{\wref}[1]{$\wlr{\ref{#1}}$}

\newcommand{\wrone}{\mathds R}
\newcommand{\wrn}[1]{{\mathds R}^{#1}}
\newcommand{\wrm}[1]{\mathrm{#1}}

\newcommand{\ws}[1]{\wcal{#1}}
\newcommand{\wset}[1]{\left\{ {#1} \right\}}
\newcommand{\wtr}[1]{\mathrm{T}}
\newcommand{\wu}[1]{\underline{#1}}
\newcommand{\wv}[1]{\mathbf{#1}}

\usepackage{authblk}
\title{Root finding with interval arithmetic}
\author[1]{Walter F. Mascarenhas}
\affil[1]{Departamento de Computação, IME\\ Universidade de São Paulo, Brazil}
\date{\vspace{-5ex}}
\begin{document}

\maketitle

\begin{abstract}
We consider the solution of nonlinear equations in one real
variable, the problem usually called by {root finding.} 
Although this is an old problem, we believe that some
aspects of its solution using interval arithmetic are
not well understood, and we present our views on
this subject. We argue that problems with just
one variable are much simpler than problems with
more variables, and we should use specific methods
for them. We provide an implementation of our ideas in \texttt{C++},
and make this code available under the 
Mozilla Public License 2.0.
\end{abstract}

\section{Introduction}
\label{sec_intro}
Many books and articles explain how to
solve the nonlinear equation $\wfc{f}{x} = 0$ for $x \in \wrone{}$,
and some of them consider the {\it verified solution} of such
equation, that is, finding solutions with rigorous bounds on them.
Here we discuss the computation of verified solutions using
interval arithmetic. This is an old subject, 
but the interval arithmetic literature
is mostly concerned with the solution of the general 
equation $\wfc{f}{x} = 0$ for $x \in \wrn{n}$ using
variations of Newton's method,
since its very beginning \cite{HansenA,RayMoore,Nickel}.  
For instance, the third chapter of
\cite{Ratschek} gives a nice description of 
interval Newton's method and \cite{HansenA,HansenB} 
present interesting ways to improve it.
Here we focus on the simplest case $n = 1$,
because it is very important in practice and there are techniques
which are applicable only for $n = 1$. As far as we know,
the most detailed discussion of this particular case
is presented in chapter 9 of \cite{HW}, but we believe
that there is more to be said about this subject
than what is presented there.  

In this article, $f$ is a function
from $\wrone{}$ to $\wrone{}$.
We denote $f$'s $\ell$th derivative by $f^{\wlr{\ell}}$ 
and when we mention such derivatives we assume that they exist. We consider the class
\[
\wds{I} := \wset{\, \wit{x} = [\wu{x},\wo{x}], \ \wrm{with} \
 -\infty \leq \wu{x} \leq \wo{x} \leq +\infty }  \cup \wset{ \emptyset }
\]
of all closed intervals,
including the empty and unbounded ones. For $\ws{S} \subset \wrone{}$,
we write
\[
\wfc{f}{\ws{S}} := \wset{\wfc{f}{x} \ \wrm{for} \ x \in \ws{S}}.
\]
We assume that we have extensions of $f$ and its derivatives, in the sense that:

\begin{wdef}
We say that a function $\wit{f}: \wds{I} \rightarrow \wds{I}$ is an extension
of a function $f: \wrone{} \rightarrow \wrone{}$ in an interval $\wit{x}$ 
if for every $\wit{x}' \subset \wit{x}$ 
we have that $\wfc{f}{\wit{x}'} \subset \wfc{\wit{f}}{\wit{x}'}$. 
\wqes{}
\end{wdef}
We usually identify the point $t \in \wrone{}$ with the interval $[t,t]$, and
for functions $\wit{f}: \wds{I} \rightarrow \wds{I}$ we write
$\wfc{\wit{f}}{t} := \wfc{\wit{f}}{[t,t]}$. The set of roots 
is denoted by $\ws{R}$, and $r$ is a generic root. 

When $f$ is differentiable, the main tool for root finding
is Newton's method
\begin{equation}
\label{newton_p}
x_{k+1} = x_k - \wfc{f}{x_k}/\wdfc{f}{x_k}.
\end{equation}
This method has a natural extension to interval arithmetic:
\begin{equation}
\label{newton_i}
\wit{x}_{k+1} = \wit{x}_k \cap \wlr{\ t_k - 
\wfc{\wit{f}}{t_k}/\, \wdfc{\wit{f}}{\wit{x}_k} \ },
\end{equation}
where $t_k$ is a point in $\wit{x}_k$ and $f$ and $f'$ are extensions of
$\wit{f}$ and $\wit{f}'$.
Traditionally, we compute $t_k$ by rounding the midpoint of $\wit{x}_k$.
The first question an alert numerical analyst would ask about Equation \wref{newton_i} is:
%%%
\begin{equation}
\label{qd0}
\wrm{What \ should \ we \ do \ when \ } \wdfc{\wit{f}}{\wit{x}_k} \ \wrm{contains} \ 0 \wrm{?}
\end{equation}
%%%
However, we see few people asking the following question:
%%%
\begin{equation}
\label{qd1}
\wrm{What \ should \ we \ do \ when \ } \wdfc{\wit{f}}{\wit{x}_k} \ 
\wrm{{\bf does \ not} \ contain} \ 0 \wrm{?}
\end{equation}
%%%

Both questions are quite relevant, but before we answer them 
we must say that the interval version of Newton's
method {\bf may be implemented without Equation \wref{newton_i}}.
Instead, we could write $f$ in the {\it centered form}
\begin{equation}
\label{center}
\wfc{f}{t} \in 
\wfc{\wit{f}}{t_k} + \wfc{\wit{s}_{t_k}}{\wit{x}_k} \wlr{\wit{x}_k - t_k}
\whs{0.5} \wrm{for} \whs{0.5} t \in \wit{x}_k,
\end{equation}
and replace the extended derivative 
$\wdfc{\wit{f}}{\wit{x}_k}$ by 
an extension $\wfc{\wit{s}_{t_k}}{\wit{x}_k}$ of the {\it slope} 
$s_{t_k}$ of  $f$ at $t_k$. This leads to an improved
version of Newton's method given by
\begin{equation}
\label{newton_t1}
\wit{x}_{k+1} = \wit{x}_k \cap \wlr{\ t_k - 
\wfc{\wit{f}}{t_k}/\,\wfc{\wit{s}_{t_k}}{\wit{x_k}} \ }.
\end{equation}
The slope is defined as 
\begin{equation}
\label{slope1}
\wfc{s_{c}}{t} := 
\wfc{s_{f,c}}{t} := 
\left\{
\begin{array}{ccc}
\wdfc{f}{c} & \wrm{if} & t = c, \\[0.1cm]
\frac{\wfc{f}{t} - \wfc{f}{c}}{t - c} & \wrm{if} & t \neq c,
\end{array}
\right.
\end{equation}
and centered forms are already
mentioned in Moore's book \cite{RayMoore}. They
have been discussed in detail in the interval 
arithmetic literature \cite{KN,RatschekC}
and have generalizations called Taylor forms
or Taylor models \cite{Neumaier}. Our algorithm
uses only the plain centered form \wref{center}, 
and in situations in which practical implementations
of such generalizations are available we would
use them only as a tool to obtain more accurate
centered forms.

%%%%%
The Mean Value Theorem shows that
\[
t \leq c \Rightarrow \wfc{s_{c}}{t} \in \wdfc{f}{[t,c]}
\whs{1} \wrm{and} \whs{1} 
t \geq c \Rightarrow
\wfc{s_{c}}{t} \in \wdfc{f}{[c,t]}
\]
and this implies that any extension of 
$\wit{f}'$ is an extension of $s_{t_k}$, but there may be better ones,
specially when the interval $\wit{x}_k$ is not small. For instance,
if $\wfc{f}{t} := t^2$ then
\[
\wfc{s_{0}}{t} = x \whs{1} \wrm{and} \whs{1}
\wdfc{f}{t} = 2 t
\]
and any extension of $f'$ yields intervals twice as
large than necessary for an slope.

In practice the gain by replacing derivatives by slopes may not be 
great, because usually $\wdfc{f}{r} \neq 0$ at the
roots $r$ and in this case, for $t_k$ close to $r$,
the difference between $s_{t_k}$ and $\wdfc{f}{t}$ 
for $t$ in a short interval $\wit{x}_k$ is not 
large. Moreover, extensions of derivatives have an important
feature that extensions of slopes do not have:
derivative extensions can detect monotonicity. 
In other words, if $\wfc{\wit{f}'}{\wit{x}} \subset (0,+\infty)$ then
$f$ is increasing in $\wit{x}$, but we cannot reach
the same conclusion from the fact that 
$\wfc{s_{t_k}}{\wit{x}} \subset (0,+\infty)$. The
ease with which we can use information about monotonicity 
when we have only one variable is what makes the 
case $n=1$ so special, and this is the reason why
this article was written to start with.

We can ask questions similar to \wref{qd0} and \wref{qd1}
about the modified Newton's step \wref{newton_t1}:
\begin{equation}
\label{qs0}
\wrm{What \ should \ we \ do \ when \ } 
\wfc{\wit{s}_{m_k}}{\wit{x}_k} \ \wrm{contains} \ 0 \wrm{?}
\end{equation}
\begin{equation}
\label{qs1}
\wrm{What \ should \ we \ do \ when \ } \wfc{\wit{s}_{m_k}}{\wit{x}_k} \ 
\wrm{does \ not \ contain} \ 0 \wrm{?}
\end{equation}
%%%
and there are at least two more questions we must ask:
\begin{equation}
\label{qret}
\wrm{What \ information \ about \ the \ roots \  our \ algorithm \ should \ provide?} 
\end{equation}
\begin{equation}
\label{qstop}
\wrm{When \ should \ we \ stop \ the \ iterations?}
\end{equation}

After much thought about the questions above, 
we devised a root finding algorithm which combines the interval
versions of Newton's method in Equations 
\wref{newton_i} and \wref{newton_t1},
and modifies them in order to exploit montonicity. Our 
algorithm tries to strike a balance between theory and
practice, and takes into account the following practical issues:
\begin{enumerate}
\item Usually, the evaluation of interval extensions over
intervals is much less accurate than the evaluation of interval extensions
at points, that is, the width of the interval $\wfc{\wit{f}}{t_k}$ 
is much smaller than the width of the interval 
$\wfc{\wit{f}}{\wit{x}_k}$. The same applies to the derivative $\wit{f}'$.
%%%
\item Usually, the floating point evaluation of $d_k = \wdfc{f}{t_k}$ yields a
reasonable estimate this derivative, at a much lower cost than
the evaluation of the extensions 
$\wdfc{\wit{f}}{t_k}$ or  $\wdfc{\wit{f}}{\wit{x}_k}$.
The only defect of $d_k$ is that 
it does not come with a guarantee of its accuracy. As a result
we can, and should, use floating point computations in order to obtain
reasonable estimates (which may turn out to be inaccurate sometimes)
and resort to interval computation only when we absolutely need guarantees
about our results. In particular, the computed $\wdfc{\wit{f}}{\wit{x}_k}$ 
may be very wide even when the floating point $d_k = \wdfc{f}{t_k}$
would lead to a good Newton step $t_{k+1} = t_k - \hat{w}_k / d_k$, where
$\hat{w}_k$ is the mid point of $\wit{w}_k = \wfc{\wit{f}}{t_k}$. 
%%%
\item In interval arithmetic, our goal is to find short intervals $\wit{r}$
which may contain roots and to discard intervals guaranteed not to contain roots.
%%%
\item
The simplest way to ensure that an interval $\wit{r} = [\wu{r},\wo{r}]$ contains a root 
is to prove that $\wfc{f}{\wu{r}}$ and $\wfc{f}{\wo{r}}$ have opposite signs.
We can do that by evaluating $\wit{f}$ at the points $\wu{r}$ and $\wo{r}$, and
as we noted above this evaluation tends do yield sharp results.
%%%
\item 
 The simplest way to ensure that an interval $\wit{r}$ does not contain a root 
is to prove that $\wfc{f}{\wu{r}}$ and $\wfc{f}{\wo{r}}$ are
different from zero and have the same sign
and $\wdfc{f}{\wit{r}}$ does not contain $0$. In practice, this is not
as easy as in the previous item because the computed $\wdfc{\wit{f}}{\wit{r}}$ tends
to be inflated by the usual weaknesses of interval arithmetic.
However, when we know that $f$ is monotonic we can dispense with 
$\wdfc{\wit{f}}{\wit{r}}$ and check only the values of $f$ at the end points. 
This simplifies things immensely for monotonic functions. Actually, it 
makes little practical sense to compute $\wdfc{\wit{f}}{\wit{x}_k}$ 
once we already know that it does not contain 0, and we have this
information for all nodes below a first node $\wit{x}_k$ in the branch and bound
tree which is such that $\wdfc{\wit{f}}{\wit{x}_k}$ does not contain 0.
%%%
\item Multiple roots, i.e., roots $r$ such that $\wdf{f}{r} = 0$, are a major
problem and should be handled differently from single roots. Usually, it is hopeless
to try to get them with the same accuracy as single roots, and the cases in which
we have some hope are too specific to deserve attention in a generic software.
\end{enumerate}
Combining the items above, we propose an algorithm which can be outlined
as follows. For each candidate interval we keep seven additional 
bits of information:
\begin{itemize}
\item The sign $\sigma_i \in \wset{-1,0,1}$ of $f$ at its infimum, 
\item The sign $\sigma_s \in \wset{-1,0,1}$ of $f$ at its supremum
\item The sign $\sigma_d \in \wset{-1,0,1}$ of $f'$ in $\wit{x}$
\item A point $t$ in its interior, indicating where it should be split in the Newton step
\item A point $\tilde{t}$, which lies near $\wit{x}$.
\item The value $\tilde{d} =\wdfc{f}{\tilde{t}}$.
\item The expected sign $\sigma_t$ of $\wfc{f}{t}$.  
\end{itemize}
Regarding the signs above, $-1$ means definitely negative, $+1$ means definitely positive
and $0$ means that we do not know the corresponding sign.

We then procedure in the usual branch and bound way:

\begin{enumerate}
\item If the stack is empty then we are done. Otherwise we pop a tuple 
\[
\wlr{\wit{x}_k,\sigma_i,\sigma_s,\sigma_d,t_k,\tilde{t}, \tilde{d}, \sigma_t}
\]
from the stack.
\item If the width of $\wit{x}$ is below a given tolerance $\tau_x$ then
we insert $\wlr{\wit{x},\sigma_i,\sigma_s, \sigma_d}$ in a list of possible solutions and go to item 1.

\item We compute $\wfc{\wit{f}}{\wit{x}}$. 
If it  does not contain $0$, then drop $\wit{x}$ and go to item 1.
(Tolerances are discussed in Section \ref{sec_expect}.)
%%%
%%%
\item We compute the slope $\wit{s}_k := \wfc{\wit{s}_{t_k}}{\wit{x}_k}$. If it does not contain
$0$ then we compute the derivative $\wit{d}_k := \wdfc{\wit{f}}{\wit{x}_k}$. If $\wit{d}_k$
does not contain $0$ then we change to the specific algorithm to handle
strictly monotonic functions described in Section \ref{sec_mono},
and once we are done we go to item 1. Otherwise, we replace $\wit{s}_k$ by 
$\wit{s}_k \cap \wit{d}_k$ and continue.
%%%
%%%
\item We compute $\wfc{\wit{f}}{t_k}$ 
and check whether it contains zero or is too close to zero. If it does
then use the algorithm described in Section \ref{sec_expect} to 
handle $\wit{x}_k$ and go to item 1.
%%%
%%%
\item If $\wfc{\wit{f}}{t_k}$  is not too close to zero then
and apply the version of the interval Newton step described in Section
\ref{sec_intern}, obtaining at most two intervals $\wit{x}_{k+1}$ and $\wit{x}_{k+1}'$. 
If there are no $\wit{x}_{k+1}$ then we drop $\wit{x}_k$ and go to item 1. 
%%%
%%%
%%%
\item If there is just one $\wit{x}_{k+1}$ then we check whether 
the sign $\sigma_t$ matches the sign of $\wfc{\wit{f}}{\wit{x}_k}$. If it does not then we set
$\sigma_t$ to zero, take $t_{k+1}$ as the mid point of $\wit{x}_{k+1}$,
obtain the extra information for $\wit{x}_{x+1}$ and push it on the stack and go to item 1.
Otherwise, we compute $d =\wdfc{f}{t_k}$ using floating point 
arithmetic and use $\tilde{d}$ and $\tilde{t}$ to check whether the corrected Newton step $t_{k+1}$ described
in Section \ref{sec_modify} lies in $\wit{x}_{k+1}$. If it does
then we use this $t_{k+1}$, otherwise we take $t_{k+1}$ as the midpoint of $\wit{x}_{k+1}$.
We then obtain the additional information for $\wit{x}_{k+1}$, push it on the stack and go to item 1.
%%%
%%%
%%%
\item 
If there are two $\wit{x}_{k+ 1}$s, then we check whether the width of 
$\wit{x}_k$ lies below the {\it cluster threshold} $\tau_c$. If it
does then we deem $\wit{x}_k$ to contain a cluster of zeros,
and proceed as in item 2.
Otherwise we apply the same procedure as in item 7 to $\wit{x}_{k+1}$ 
and $\wit{x}_{k+1}'$  and go to item 1.
\end{enumerate}

We implemented the algorithm outlined above in \texttt{C++},
using our Moore interval arithmetic library \cite{Moore}.
This code is available with
the arxiv version of this article. It is distributed under the
Mozilla Public License 2.0. Unfortunately, the code is
much more involved than the summary above, because
it must deal with many technical details which
were not mentioned in this article in order not to
make it longer and more complex than it already is.

In the rest of the article we discuss in more detail
several aspects of the algorithm outline above.
We start with Section \ref{sec_intern}, in which
we present a version of Newton's method for interval arithmetic.
This version is similar to the ones found in the literature,
but it is slightly different because it ensures that
$f$ is different from zero at some of the extreme points of the new
intervals, and computes the signs of $f$ at these extremes at a low extra cost. 
As a result, our algorithm yields not only intervals
containing the roots but also the sign of $f$ at the extreme of such intervals,
and these signs certify the existence of a root in the interval when they are different
from each other. In Section \ref{sec_expect}  we make some comments about
interval arithmetic in general and discuss the thorny subject of tolerances,
which are unavoidable for defining stopping criterion for interval root solvers.
Section \ref{sec_modify} presents yet
another version of the classic Newton's method
for exact, ``point arithmetic.'' This ``point version''
is the motivation for the interval version of Newton's 
method for monotone functions presented in Section
\ref{sec_mono}. Finally, in Section \ref{sec_test} we discuss the important
subject of testing.

%%%%%%%%%%%%%
%%%%%%%%%%%%%
%%%%%%%%%%%%%
%%%%%%%%%%%%%

\section{The Interval version of Newton's step}
\label{sec_intern}
This section is about the interval version of Newton's step
\begin{equation}
\label{newton_t2}
\wit{x}_{k+1} = \wit{x}_k \cap \wlr{\ t_k - 
\wfc{\wit{f}}{t_k}/\, \wit{d}_k \ },
\end{equation}
where $\wit{d}_k = [\wu{d}_k, \wo{d}_k]$ can be either the derivative $\wdfc{\wit{f}}{\wit{x}_k}$ 
or the slope $\wfc{\wit{s}_{t_k}}{\wit{x_k}}$. Here we make the simplifying
assumption that the interval 
\[
\wit{w}_k := \wfc{\wit{f}}{t_k} := [\wu{w}_k, \ \wo{w}_k]  
\]
does not contain $0$, answering to the questions 
\wref{qd0}, \wref{qd1}, \wref{qs0} and \wref{qs1} in the introduction in this case.
By replacing $f$ by $-f$ if necessary, we can assume that $\wo{w}_k < 0$,
and we make this simplifying assumption from now on.

The answer to questions \wref{qd0} and \wref{qs0}, in which case
$\wit{d}_k$ contains $0$, is illustrated in Figure \ref{figsw0}.
In this Figure, the inclined lines have equations
\begin{equation}
\label{inter1}
w = \wo{w}_k + \wu{d}_k \wlr{t - t_k} \whs{1} \wrm{and} \whs{1}
w = \wo{w}_k + \wo{d}_k \wlr{t - t_k},
\end{equation}
and by intersecting these lines with the axis $w = 0$
we obtain a better bound on $\ws{R} \cap \wit{x}_k$ 
(recall that $\ws{R}$ is the set of roots.)
There are three possibilities for our bounds on $\ws{R} \cap \wit{x}_k$ after
we take in to account the intersections above. We may find that:
%%%%%%
\begin{itemize}
\item $\ws{R} \cap \wit{x}_k = \emptyset$. In this case we drop $\wit{x}_k$.
%%%%%%
\item $\ws{R} \cap \wit{x}_k$ is contained in a single 
interval $\wit{r} = [\wu{r},\wo{r}]$, as in the second and third cases in 
Figure \ref{figw0}. In this case we take $\wit{x}_{k+1} = \wit{r}$.
%%%%%%
\item $\ws{R} \cap \wit{x}_k $ is contained in the union of two 
intervals $\wit{r}_1$ and $\wit{r}_2$ such that $\wo{r}_1 \leq \wu{r}_2$,
as in the last case of Figure \ref{figw0}. 
\end{itemize}

\begin{figure}[!h]
\begin{tikzpicture}[scale = 1]

%%%%%%%%%%%%
%% small low, small up
%%%%%%%%%%%%

\draw node at (1.3, 2.2){$\ws{R} \cap \wit{x}_k$ is empty};
\draw (0.0, 2.0)--(2.8, 2.0);

\filldraw[black] (0.0, 2.0) circle(1.5pt);
\filldraw[black] (2.8, 2.0) circle(1.5pt);

\draw (0.1, 2.2)--(0.0, 2.2)--( 0.0, 1.8)--(0.1, 1.8);
\draw (2.7, 2.2)--(2.8, 2.2)--( 2.8, 1.8)--(2.7, 1.8);

\draw (1.4, 0.0)--(1.4, 1.0);
\filldraw[black] (1.4, 0.0) circle(1.5pt);
\filldraw[black] (1.4, 1.0) circle(1.5pt);
\draw (1.2, 0.1)--(1.2, 0.0)--(1.6, 0.0)--(1.6, 0.1);
\draw (1.2, 0.9)--(1.2, 1.0)--(1.6, 1.0)--(1.6, 0.9);
\draw node at (1.45, -0.3){$t_k$};

\draw node at (1.9, 0.9){$\wo{w}_k$};
\draw node at (1.9, 0.1){$\wu{w}_k$};

\draw (1.4, 1.0)--(0.0, 1.5);
\draw (1.4, 1.0)--(2.8, 1.5);

%%%%%%%%%%%%
%% small low, large up
%%%%%%%%%%%%

\draw (3.1, 2.0)--( 5.9, 2.0);

\filldraw[black] (3.1, 2.0) circle(1.5pt);
\filldraw[black] (5.9, 2.0) circle(1.5pt);

\draw (3.2, 2.2)--(3.1, 2.2)--( 3.1, 1.8)--(3.2, 1.8);
\draw (5.8, 2.2)--(5.9, 2.2)--( 5.9, 1.8)--(5.8, 1.8);

\draw (4.5, 0.0)--(4.5, 1.0);
\filldraw[black] (4.5, 0.0) circle(1.5pt);
\filldraw[black] (4.5, 1.0) circle(1.5pt);
\draw (4.3, 0.1)--(4.3, 0.0)--(4.7, 0.0)--(4.7, 0.1);
\draw (4.3, 0.9)--(4.3, 1.0)--(4.7, 1.0)--(4.7, 0.9);
\draw node at (4.5, -0.3){$t_k$};

\draw node at (5.0, 0.9){$\wo{w}_k$};
\draw node at (5.0, 0.1){$\wu{w}_k$};

\filldraw[black] (5.25, 2.0) circle(1.5pt);
\draw node at (5.25, 1.6){$\wu{r}$};
\draw node at (5.9,  1.6){$\wo{r}$};

\draw (4.5, 1.0)--(3.1, 1.5);
\draw (4.5, 1.0)--(5.9, 2.9);

%%%%%%%%%%%%
%% large up low, small up
%%%%%%%%%%%%

\draw (6.2, 2.0)--(9.0, 2.0);

\filldraw[black] (6.2, 2.0) circle(1.5pt);
\filldraw[black] (9.0, 2.0) circle(1.5pt);

\draw (6.3, 2.2)--(6.2, 2.2)--( 6.2, 1.8)--(6.3, 1.8);
\draw (8.9, 2.2)--(9.0, 2.2)--( 9.0, 1.8)--(8.9, 1.8);

\draw (7.6, 0.0)--(7.6, 1.0);
\filldraw[black] (7.6, 0.0) circle(1.5pt);
\filldraw[black] (7.6, 1.0) circle(1.5pt);
\draw (7.4, 0.1)--(7.4, 0.0)--(7.8, 0.0)--(7.8, 0.1);
\draw (7.4, 0.9)--(7.4, 1.0)--(7.8, 1.0)--(7.8, 0.9);
\draw node at (7.6, -0.3){$t_k$};
\draw node at (8.1, 0.9){$\wo{w}_k$};
\draw node at (8.1, 0.1){$\wu{w}_k$};

\draw (7.6, 1.0)--(6.2, 2.9);
\draw (7.6, 1.0)--(9.0, 1.5);

\draw node at (6.2,  2.45){$\wu{r}$};
\draw node at (6.85, 2.45){$\wo{r}$};
\filldraw[black] (6.85, 2.0) circle(1.5pt);

%%%%%%%%%%%%
%% large low , large up
%%%%%%%%%%%%

\draw (9.3, 2.0)--(12.1, 2.0);

\filldraw[black] ( 9.3, 2.0) circle(1.5pt);
\filldraw[black] (12.1, 2.0) circle(1.5pt);
\draw ( 9.4, 2.2)--( 9.3, 2.2)--(  9.3, 1.8)--( 9.4, 1.8);
\draw (12.0, 2.2)--(12.1, 2.2)--( 12.1, 1.8)--(12.0, 1.8);

\draw (10.7, 0.0)--(10.7, 1.0);
\filldraw[black] (10.7, 0.0) circle(1.5pt);
\filldraw[black] (10.7, 1.0) circle(1.5pt);
\draw (10.5, 0.1)--(10.5, 0.0)--(10.9, 0.0)--(10.9, 0.1);
\draw (10.5, 0.9)--(10.5, 1.0)--(10.9, 1.0)--(10.9, 0.9);

\draw node at (11.2, 0.9){$\wo{w}_k$};
\draw node at (11.2, 0.1){$\wu{w}_k$};
\draw node at (10.7, -0.3){$t_k$};

\draw (10.7, 1.0)--( 9.3, 2.9);
\draw (10.7, 1.0)--(12.1, 2.9);

\draw node at (9.3,  2.45){$\wu{r}_1$};
\draw node at (9.95, 2.45){$\wo{r}_1$};
\filldraw[black] (9.95, 2.0) circle(1.5pt);

\filldraw[black] (11.45, 2.0) circle(1.5pt);
\draw node at (11.5, 1.6){$\wu{r}_2$};
\draw node at (12.1,  1.6){$\wo{r}_2$};

\end{tikzpicture}
\caption{The case $0 \in \wit{d}$, $\wo{w} < 0$, with $-\wu{d}_k$ and $\wo{d}_k$ large and small.} 
\label{figw0}
\end{figure}

There is nothing new in our approach up to this point. We just
explained with a picture what most people do algebraically.
Our slight improvement comes in the way we compute the intersections
in Equations \wref{inter1}.  These intersections are
\[
\wo{r}_1 := t_k - \wo{w}_k / \wu{d}_k 
\hspace{1cm} \wrm{and} \hspace{1cm}
\wu{r}_2 := t_k - \wo{w}_k / \wo{d}_k,
\]
and, as usual in the Moore Library \cite{Moore}, we
propose that $\wo{r}_1$ and $\wu{r}_2$ be computed with the rounding
mode upwards, using the following expressions:
\begin{eqnarray}
\label{computedr1}
\wo{r}_1 & := &t_k + \wo{q}_k \whs{0.5} \wrm{for} \whs{0.3} \wo{q}_
k := \wo{w}_k / \wlr{-\wu{d}_k}, \\
\label{computedr2}
\wu{r}_2 & := & - u_k   
\whs{1.0} \wrm{for} \whs{0.3} u_k := \wu{q}_k - t_k \whs{0.3} \wrm{and} \whs{0.3}
\wu{q}_k := \wo{w}_k / \wo{d}_k.
\end{eqnarray}
This order of evaluation ensure that $\wo{r}_1$ and $\wu{r}_2$ are
rounded in the correct direction, so that we do not risk loosing roots.
The expressions in Equation \wref{computedr1} and \wref{computedr2} are so simple that we
find whether $\wo{r}_1$ and $\wu{r}_2$ where computed exactly by checking
whether
\begin{equation}
\label{veri1}
\wo{r}_1 -  t_k = \wo{q}_k\whs{0.5} \wrm{and} \whs{0.5} \wo{q}_k \wu{d}_k = \wo{w}_k
\end{equation}
and
\begin{equation}
\label{veri2}
u_k + t_k = \wu{q}_k \whs{0.5} \wrm{and} \whs{0.5} \wu{q}_k \wo{d}_k = \wo{w}_k.
\end{equation}
If we find that $\wo{r}_1$ was computed exactly then we increase it to the next floating point number. If $\wu{r}_2$ was computed exactly then we decrease it to the previous floating point
number. By doing so, we ensure that $\wfc{f}{\wo{r}_1} < 0$ and $\wfc{f}{\wu{r}_2} < 0$,
without computing $\wfc{\wit{f}}{\wo{r}_1}$ or $\wfc{\wit{f}}{\wu{r}_2}$.
We should also mention that
it is possible to prove that, even after the rounding, incrementing and decrementing
steps above, $\wo{r}_1 \leq t_k \leq \wu{r}_2$ and there is no risk of 
$\wo{r}_1$ crossing $\wu{r}_2$.

Regarding the cost of all this, note that in the Moore library the rounding mode is
always upwards. Therefore, there is no cost associated with changing rounding modes.
Moreover, the expensive operations in Equations \wref{computedr1} and \wref{computedr2}
are the divisions,
and the extra couple of sums and multiplications do not increase the overall cost
of computing the intersections by much. In fact, if we take into account branch
prediction and speculative execution 
and the fact that computations are usually inexact, the extra cost
due to the verification in Equations \wref{veri1} and \wref{veri2} and
the decrement of $\wo{r}_1$ and the increment of $\wu{r}_2$ is likely
to be minimal. We would not be surprised if by performing the
verification in Equations \wref{veri1} and \wref{veri2}
above the code would be faster than blindly incrementing $\wo{r}_1$
and decrementing $\wu{r}_2$ (we did not have the time to check this.)

Finally, by using symmetry, we can reduce 
the analysis of questions \wref{qd1} and \wref{qs1}
in the case in which $\wfc{\wit{f}}{t_k}$ does not contain $0$
to the cases described in Figure \ref{figsw0}.
In this Figure, the inclined lines have equations
\[
w = \wo{w}_k + \wo{d}_k \wlr{t - t_k} \whs{1} \wrm{and} \whs{1}
w = \wu{w}_k + \wu{d}_k \wlr{t - t_k}.
\] 
and we may either have no intersection or one intersection. We can
use the same technique described above for find intersections
which are correctly rounded and such that $\wfc{f}{u}$ is
different from zero in the new extreme points $u$,
without evaluating $\wfc{\wit{f}}{u}$  (the old
extreme points stay as they were.)

\begin{figure}[!h]
\begin{tikzpicture}[scale = 1]

%%%%%%%%%%%%
%% small low, small up
%%%%%%%%%%%%

\draw node at (1.4, 2.2){$\ws{R} \cap \wit{x}_k$ is empty};
\draw (0.0, 2.0)--(2.8, 2.0);

\filldraw[black] (-1.5, 2.0) circle(0.01pt);

\filldraw[black] (0.0, 2.0) circle(1.5pt);
\filldraw[black] (2.8, 2.0) circle(1.5pt);

\draw (0.1, 2.2)--(0.0, 2.2)--( 0.0, 1.8)--(0.1, 1.8);
\draw (2.7, 2.2)--(2.8, 2.2)--( 2.8, 1.8)--(2.7, 1.8);

\draw (1.4, 0.0)--(1.4, 1.0);
\filldraw[black] (1.4, 0.0) circle(1.5pt);
\filldraw[black] (1.4, 1.0) circle(1.5pt);
\draw (1.2, 0.1)--(1.2, 0.0)--(1.6, 0.0)--(1.6, 0.1);
\draw (1.2, 0.9)--(1.2, 1.0)--(1.6, 1.0)--(1.6, 0.9);
\draw node at (1.45, -0.3){$t_k$};

\draw node at (0.9, 0.9){$\wo{w}_k$};
\draw node at (0.9, 0.1){$\wu{w}_k$};

\draw (1.4, 0.0)--(2.8, 0.5);
\draw (1.4, 1.0)--(2.8, 1.5);

%%%%%%%%%%%%
%% small low, large up
%%%%%%%%%%%%

\draw (3.1, 2.0)--( 5.9, 2.0);

\filldraw[black] (3.1, 2.0) circle(1.5pt);
\filldraw[black] (5.9, 2.0) circle(1.5pt);

\draw (3.2, 2.2)--(3.1, 2.2)--( 3.1, 1.8)--(3.2, 1.8);
\draw (5.8, 2.2)--(5.9, 2.2)--( 5.9, 1.8)--(5.8, 1.8);

\draw (4.5, 0.0)--(4.5, 1.0);
\filldraw[black] (4.5, 0.0) circle(1.5pt);
\filldraw[black] (4.5, 1.0) circle(1.5pt);
\draw (4.3, 0.1)--(4.3, 0.0)--(4.7, 0.0)--(4.7, 0.1);
\draw (4.3, 0.9)--(4.3, 1.0)--(4.7, 1.0)--(4.7, 0.9);
\draw node at (4.5, -0.3){$t_k$};

\draw node at (4.0, 0.9){$\wo{w}_k$};
\draw node at (4.0, 0.1){$\wu{w}_k$};

\filldraw[black] (4.9, 2.0) circle(1.5pt);
\draw node at (5.0, 1.6){$\wu{r}$};
\draw node at (5.9,  1.6){$\wo{r}$};

\draw (4.5, 0.0)--(5.9, 1.0);
\draw (4.5, 1.0)--(5.3, 2.9);

%%%%%%%%%%%%
%% large up low, small up
%%%%%%%%%%%%

\draw (6.2, 2.0)--(9.0, 2.0);

\filldraw[black] (6.2, 2.0) circle(1.5pt);
\filldraw[black] (9.0, 2.0) circle(1.5pt);

\draw (6.3, 2.2)--(6.2, 2.2)--( 6.2, 1.8)--(6.3, 1.8);
\draw (8.9, 2.2)--(9.0, 2.2)--( 9.0, 1.8)--(8.9, 1.8);

\draw (7.6, 0.0)--(7.6, 1.0);
\filldraw[black] (7.6, 0.0) circle(1.5pt);
\filldraw[black] (7.6, 1.0) circle(1.5pt);
\draw (7.4, 0.1)--(7.4, 0.0)--(7.8, 0.0)--(7.8, 0.1);
\draw (7.4, 0.9)--(7.4, 1.0)--(7.8, 1.0)--(7.8, 0.9);
\draw node at (7.6, -0.3){$t_k$};
\draw node at (7.1, 0.9){$\wo{w}_k$};
\draw node at (7.1, 0.1){$\wu{w}_k$};

\filldraw[black] (7.9, 2.0) circle(1.5pt);
\draw (7.6, 1.0)--(8.2, 2.9);
\draw (7.6, 0.0)--(9.0, 2.9);
\filldraw[black] (8.55, 2.0) circle(1.5pt);

\draw node at (7.8, 2.35){$\wu{r}$};
\draw node at (8.5, 2.35){$\wo{r}$};

\end{tikzpicture}
\caption{The case $\wu{d}_k > 0$ and $\wo{w}_k < 0$, with $\wu{d}_k$ and $\wo{d}_k$ large and small.} 
\label{figsw0}
\end{figure}

\section{What we can expect from Interval Arithmetic}
\label{sec_expect}
In order to appreciate the content of this article,
one must realize that interval arithmetic is quite different from 
floating point arithmetic, and it is used for different purposes.
Floating point methods are usually faster than interval methods, but
do not provide rigorous bounds on their results. We use
interval arithmetic when we want global solutions to our
problems, with rigorous bounds on them, and we are willing
to pay more for that. 

\begin{figure}[!h]
\begin{tikzpicture}[scale = 1]

\draw (0,0)--(12,0);
\draw (0.1,0.1)--(0,0.1)--(0,-0.1)--(0.1,-0.1);
\filldraw[black] (0,0) circle(1pt);
\filldraw[black] (12,0) circle(1pt);
\draw (11.9,0.1)--(12,0.1)--(12,-0.1)--(11.9,-0.1);

\draw (0,-0.6) .. controls (6,-0.25) .. (12,0.5);
\filldraw[black] (0, -0.6) circle(1pt);
\filldraw[black] (12, 0.5) circle(1pt);

\draw (4.86,-0.1 )--(4.76,-0.1)--(4.76,0.1)--(4.86,0.1);
\draw node at (4.76,0.3){$\wu{r}$};
\filldraw[black](4.76,0) circle(1pt);
\filldraw[black](7.93,0) circle(1pt);
\filldraw[black](10.67,0) circle(1pt);
\draw node at (10.67,-0.3){$\wo{r}$};
\draw (10.57,-0.1 )--(10.67,-0.1)--(10.67,0.1)--(10.57,0.1);

\draw (5.9,-0.1 )--(5.8,-0.1)--(5.8,0.1)--(5.9,0.1);
\draw node at (5.8,0.45){$\wu{r}'$};

\draw (5.8, 0)--(5.8, 0.2)--(7, 0.2)--(7, 0.4)--(9.2, 0.4)--(9.2,0.6)--(11.5,0.6)--(11.5,0.8)--(12,0.8);
\filldraw[black] (12,0.8) circle(1pt);
\filldraw[black] (11.5,0.8) circle(1pt);
\draw[black] (11.5,0.6) circle(1pt);
\draw[black] (9.2,0.6) circle(1pt);
\filldraw[black] (9.2,0.4) circle(1pt);

\filldraw(7, 0.4) circle(1pt);
\draw(7, 0.2) circle(1pt);
\draw(5.8, 0.2) circle(1pt);
\filldraw(5.8, 0) circle(1pt);
\draw (4.76,0)-- (4.76,-0.2)--(2.43,-0.2)--(2.43, -0.4)--(0,-0.4);
\draw (4.76, -0.2) circle(1pt);
\filldraw[black] (2.43, -0.2) circle(1pt);
\draw[black] (2.43, -0.4) circle(1pt);
\filldraw[black] (0, -0.4) circle(1pt);

\draw (0,-0.8)--(3.8,-0.8)--(3.8, -0.6)--(5.5, -0.6)--(5.5,-0.4)--(8.7,-0.4)--(8.7,-0.2)--(9.7,-0.2)--(9.7,0)--(10.67,0)--(10.67,0)--(10.67,0.2)--(12,0.2);

\draw (9.6,-0.1 )--(9.7,-0.1)--(9.7,0.1)--(9.6,0.1);
\draw node at (9.7,-0.45){$\wo{r}'$};

\filldraw(0,  -0.8) circle(1pt);
\filldraw(3.8,-0.8) circle(1pt);
\draw(    3.8,-0.6) circle(1pt);
\filldraw(5.5,-0.6) circle(1pt);
\draw(    5.5,-0.4) circle(1pt);
\filldraw(8.7, -0.4) circle(1pt);
\draw(8.7, -0.2) circle(1pt);
\draw(9.7, -0.2) circle(1pt);
\filldraw(9.7,  0) circle(1pt);
\draw(10.67, 0.2) circle(1pt);
\draw(12,  0.2) circle(1pt);

\draw (7.93,-0.4)-- (7.93,0.4);
\filldraw[black](7.93,-0.4) circle(1pt);
\filldraw[black](7.93, 0.4) circle(1pt);
\draw (7.83,-0.3 )--(7.83,-0.4)--(8.03,-0.4)--(8.03,-0.3);
\draw (7.83, 0.3)--(7.83,0.4)--(8.03,0.4)--(8.03,0.3);

\draw node at ( 9.8,0.4){$\wfc{f}{t}$};
\draw node at ( 6.2,-0.7){$\wfc{\wu{\wit{f}}}{t}$};
\draw node at ( 8.6, 0.7){$\wfc{\wo{\wit{f}}}{t}$};
\draw node at (7.93,-0.6){$r$};

\end{tikzpicture}
\caption{If $\wit{f}$ is all we know, then the best estimate we can hope for $r$ is the interval 
$\wit{r} = [\wu{r}, \wo{r}]$. Knowing that 
$\wdfc{\wu{\wit{f}}}{t} > 0$ for all $t \in \wit{r}$
we can improve this to $\wit{r}' = [\wu{r}',\wo{r}']$.
This is the best we can do if $\wit{f}$ an $\wit{f}'$ are
all the information we have about $f$.}
\label{figg}
\end{figure}

Interval arithmetic is fundamentally limited.
When we evaluate the extension $\wit{f}$ at a point $t$
obtain an interval 
$\wit{w} = \wfc{\wit{f}}{t}$  containing $\wfc{f}{t}$. 
The width of $\wit{w}$ depends on the quality of the implementation of $\wit{f}$
and the underlying arithmetic. When $f$ is not too complex
and $\wfc{f}{t}$ is $\wfc{O}{1}$, we expect that the width of 
$\wit{w}$ to be $\wfc{O}{\epsilon}$,
where $\epsilon$ is the machine precision. It
is also reasonable to expect that the functions
$\wu{\wit{f}},\wo{\wit{f}}: \wit{x} \rightarrow \wrone{}$ 
will be as in Figure \ref{figg} when 
$r \in \wit{x}$ is a root of $f$.

When the evaluation of $f$ using interval arithmetic 
leads to wide intervals, we can (and do) use centered forms and Taylor models
to reduce the width of $\wfc{\wit{f}}{\wit{x}}$, but 
this will not free us from facing the reality that the evaluation
of interval extensions is usually imperfect. 
This is a fact, and we need to be prepared to handle not only
the monotone situation described in Figure \ref{figg}, but also
the noisy situation described in Figure \ref{figc}, which
often occurs near a double root, and gets worse near roots
of higher multiplicities.

\begin{figure}[!h]
\begin{tikzpicture}[scale = 1]

\draw node at ( 0.6,4.0){$\wfc{f_1}{t}$};
\draw node at (11.4,4.0){$\wfc{f_2}{t}$};
\draw node at (11.7,2.3){$\wfc{\wu{\wit{f}}}{t}$};
\draw node at ( 6.0,2.4){$\wfc{\wo{\wit{f}}}{t}$};

\draw (0,0.5)--(12,0.5);
\draw (0,0.5)--(12,0.5);
\draw [dashed] (0,1.7) -- (12,1.7);

\draw node at (0.2,1.1){$\tau_w$};
\draw[->]     (0.2,1.3)-- (0.2,1.7);
\draw[->]     (0.2,0.9)-- (0.2,0.5);

\draw (0.0, 3.5)--
			(1 * 0.15, 3.5 - 1 * 0.18 - 0.2)--(2 * 0.15, 3.5 - 2 * 0.18 + 0.2)--
			(3 * 0.15, 3.5 - 3 * 0.18 - 0.2)--(4 * 0.15, 3.5 - 4 * 0.18 + 0.2)--
			(5 * 0.15, 3.5 - 5 * 0.18 - 0.2)--(6 * 0.15, 3.5 - 6 * 0.18 + 0.2)--
	  	(7 * 0.15, 3.5 - 7 * 0.18 - 0.2)--(8 * 0.15, 3.5 - 8 * 0.18 + 0.2)--
	  	(9 * 0.15, 3.5 - 9 * 0.18 - 0.2)--
      (1.5, 1.7)--
      (1.5 + 1 * 0.12, 1.7 - 1 * 0.09 - 0.2)--(1.5 + 2 * 0.12, 1.7 - 2 * 0.09 + 0.18)--
      (1.5 + 3 * 0.12, 1.7 - 3 * 0.09 - 0.2)--(1.5 + 4 * 0.12, 1.7 - 4 * 0.09 + 0.18)--
      (1.5 + 5 * 0.12, 1.7 - 5 * 0.09 - 0.2)--(1.5 + 6 * 0.12, 1.7 - 6 * 0.09 + 0.18)--            
      (1.5 + 7 * 0.12, 1.7 - 7 * 0.09 - 0.2)--(1.5 + 8 * 0.12, 1.7 - 8 * 0.09 + 0.18)--            
      (1.5 + 9 * 0.12, 1.7 - 9 * 0.09 - 0.2)--
      (2.7, 0.8)--
      (2.7 + 1 * 0.13, 0.8 - 1 * 0.07 - 0.2)--
      (2.7 + 2 * 0.13, 0.8 - 2 * 0.07 - 0.2)--(2.7 + 3 * 0.13, 0.8 - 3 * 0.07 + 0.16)--
      (2.7 + 4 * 0.13, 0.8 - 4 * 0.07 - 0.2)--(2.7 + 5 * 0.13, 0.8 - 5 * 0.07 + 0.16)--
      (2.7 + 6 * 0.13, 0.8 - 6 * 0.07 - 0.2)--(2.7 + 7 * 0.13, 0.8 - 7 * 0.07 + 0.16)--
      (2.7 + 8 * 0.13, 0.8 - 8 * 0.07 - 0.2)--(2.7 + 9 * 0.13, 0.8 - 9 * 0.07 + 0.16)--
      (4.0, 0.1)--
      (6.0 - 9 * 0.2, 0.6 - 0.05 - 0.2)--(6.0 - 8 * 0.2, 0.6 - 2 * 0.05 + 0.0)--
      (6.0 - 7 * 0.2, 0.6 - 0.05 - 0.2)--(6.0 - 6 * 0.2, 0.6 - 2 * 0.05 + 0.1)--
      (6.0 - 5 * 0.2, 0.6 - 0.05 - 0.2)--(6.0 - 4 * 0.2, 0.6 - 2 * 0.05 + 0.1)--
      (6.0 - 3 * 0.2, 0.6 - 0.05 - 0.2)--(6.0 - 2 * 0.2, 0.6 - 2 * 0.05 + 0.1)--
      (6.0 - 1 * 0.2, 0.6 - 0.05 - 0.2)--
      ( 6.0, 0.6)--
      (6.0 + 1 * 0.2, 0.6 - 0.05 - 0.2)--(6.0 + 2 * 0.2, 0.6 - 2 * 0.05 + 0.1)--
      (6.0 + 3 * 0.2, 0.6 - 0.05 - 0.2)--(6.0 + 4 * 0.2, 0.6 - 2 * 0.05 + 0.1)--
      (6.0 + 5 * 0.2, 0.6 - 0.05 - 0.2)--(6.0 + 6 * 0.2, 0.6 - 2 * 0.05 + 0.1)--
      (6.0 + 7 * 0.2, 0.6 - 0.05 - 0.2)--(6.0 + 8 * 0.2, 0.6 - 2 * 0.05 + 0.0)--
      (6.0 + 9 * 0.2, 0.6 - 0.05 - 0.2)--
      ( 8.0, 0.1)--
      ( 8.0 + 1 * 0.13, 0.1 + 1 * 0.07 - 0.2)--( 8.0 + 2 * 0.13, 0.1 + 2 * 0.07 + 0.16)-- 
      ( 8.0 + 3 * 0.13, 0.1 + 3 * 0.07 - 0.2)--( 8.0 + 4 * 0.13, 0.1 + 4 * 0.07 + 0.16)-- 
      ( 8.0 + 5 * 0.13, 0.1 + 5 * 0.07 - 0.2)--( 8.0 + 6 * 0.13, 0.1 + 6 * 0.07 + 0.16)-- 
      ( 8.0 + 7 * 0.13, 0.1 + 7 * 0.07 - 0.2)--( 8.0 + 8 * 0.13, 0.1 + 8 * 0.07 + 0.16)-- 
      ( 8.0 + 9 * 0.13, 0.1 + 9 * 0.07 - 0.2)--
      ( 9.3, 0.8)--
      ( 9.3 + 1 * 0.12, 0.8 + 1 * 0.09 - 0.2)--( 9.3 + 2 * 0.12, 0.8 + 2 * 0.09 + 0.2)-- 
      ( 9.3 + 3 * 0.12, 0.8 + 3 * 0.09 - 0.2)--( 9.3 + 4 * 0.12, 0.8 + 4 * 0.09 + 0.2)-- 
      ( 9.3 + 5 * 0.12, 0.8 + 5 * 0.09 - 0.2)--( 9.3 + 6 * 0.12, 0.8 + 6 * 0.09 + 0.2)-- 
      ( 9.3 + 7 * 0.12, 0.8 + 7 * 0.09 - 0.2)--( 9.3 + 8 * 0.12, 0.8 + 8 * 0.09 + 0.2)-- 
      ( 9.3 + 9 * 0.12, 0.8 + 9 * 0.09 - 0.2)--
      (10.5, 1.7)--
      (10.5 + 1 * 0.15, 1.7 + 1 * 0.18 - 0.2)--(10.5 + 2 * 0.15, 1.7 + 2 * 0.18 + 0.2)--
      (10.5 + 3 * 0.15, 1.7 + 3 * 0.18 - 0.2)--(10.5 + 4 * 0.15, 1.7 + 4 * 0.18 + 0.2)--
      (10.5 + 5 * 0.15, 1.7 + 5 * 0.18 - 0.2)--(10.5 + 6 * 0.15, 1.7 + 6 * 0.18 + 0.2)--
      (10.5 + 7 * 0.15, 1.7 + 7 * 0.18 - 0.2)--(10.5 + 8 * 0.15, 1.7 + 8 * 0.18 + 0.2)--
      (10.5 + 9 * 0.15, 1.7 +  9 * 0.18 - 0.2)--
      (12.0, 3.5);

\draw (    0.0, 5.5)--
      (1 * 0.6, 5.5 - 1 * 0.35 - 0.2)-- (2 * 0.6, 5.5 - 2 * 0.35 + 0.2)--
      (3 * 0.6, 5.5 - 3 * 0.35 - 0.2)-- (4 * 0.6, 5.5 - 4 * 0.35 + 0.2)--
      (5 * 0.6, 5.5 - 5 * 0.35 - 0.2)-- (6 * 0.6, 5.5 - 6 * 0.35 + 0.2)--
      (7 * 0.6, 5.5 - 7 * 0.35 - 0.2)-- (8 * 0.6, 5.5 - 8 * 0.35 + 0.2)--
      (9 * 0.6, 5.5 - 9 * 0.35 - 0.2)-- 
      (6.0          , 2.0)-- 
      (6.0 + 1 * 0.6, 2.0 + 1 * 0.35 - 0.2)-- 
      (6.0 + 2 * 0.6, 2.0 + 2 * 0.35 + 0.2)--(6.0 + 3 * 0.6, 2.0 + 3 * 0.35 - 0.2)-- 
      (6.0 + 4 * 0.6, 2.0 + 4 * 0.35 + 0.2)--(6.0 + 5 * 0.6, 2.0 + 5 * 0.35 - 0.2)-- 
      (6.0 + 6 * 0.6, 2.0 + 6 * 0.35 + 0.2)--(6.0 + 7 * 0.6, 2.0 + 7 * 0.35 - 0.2)-- 
      (6.0 + 8 * 0.6, 2.0 + 8 * 0.35 + 0.2)--(6.0 + 9 * 0.6, 2.0 + 9 * 0.35 - 0.2)-- 
			(12.0, 5.5);

\draw node at (0.2,-0.4){$\wu{x}$};
\draw node at (1.5,-0.4){$\wu{c}$};
\draw node at (3.0,0.0){$\tau_c$};
\draw[->]     (2.7,0.0)-- (1.5,0.0);
\draw[->]     (3.2,0.0)-- (4.5,0.0);
\draw[dashed] (1.5,0.0)-- (1.5,5.5);

\draw node at (4.5,-0.4){$z$};
\draw node at (6.0,0.0){$\tau_c$};
\draw[->]     (5.7,0.0)-- (4.5,0.0);
\draw[->]     (6.2,0.0)-- (7.5,0.0);
\draw[dashed] (4.5,0.0)-- (4.5,5.5);

\draw node at (9.0,0.0){$\tau_c$};
\draw[->]     (8.7,0.0)-- (7.5,0.0);
\draw[->]     (9.2,0.0)-- (10.5,0.0);
\draw[dashed] ( 7.5,0.0)-- ( 7.5, 5.5);
\draw[dashed] (10.5,0.0)-- (10.5, 5.5);
\draw node at (10.5,-0.4){$\wo{c}$};
\draw node at (11.9,-0.4){$\wo{x}$};

\draw (0,4) .. controls (4.5,-1.5) and (3.5, 0.0) .. (12,5);
\draw (0,5) .. controls (8.5, 0.0) and (7.5,-1.5) .. (12,4);

\end{tikzpicture}
\caption{The hardest case our algorithm faces: a cluster of zeros
below the accuracy of $\wit{f}$ and $\wit{f}'$. There are no silver
bullets in such cases, and we
must find reasonable ways to cope with the $\wit{f}$ and $\wit{f}'$
which are provided to us, not the ones that we dream about. 
Using the tolerances $\tau_c$ and $\tau_w$ described in the text
is a reasonable compromise to handle such situations (which is far
from perfect.)}
\label{figc}
\end{figure}

Another fundamental fact about interval arithmetic is that it is a
pessimistic theory: it always considers the worst scenario,
and it would be inconsistent if it were otherwise.
As consequence, of this fundamental fact,
the functions $f_1$ and $f_2$ are indistinguishable by
their common extension $\wit{f}$, that is, if $\wit{f}$
is all the information that we have, then all our
conclusions about $f_1$ must also apply to $f_2$. 
This forces us to be very conservative, and include
all possible roots for $f_1$ as candidates for roots
for $f_2$ too, and vice versa. We emphasize this point
because we found that it is not always clear in the minds of 
people trying to implement roots finders like
the one we propose here, and this makes them
under estimate the effort required for this task in real life.

Taking into account all that was said since the beginning of this
section, we reached the conclusion that we should use three
tolerances to try to handle the problems caused by the very
nature of interval arithmetic. First, we should have a tolerance
$\tau_x$ to handle ``normal'' roots $r$ like the one
in Figure \ref{figg}, at which $\wdfc{f}{r} \neq 0$, so
that we consider intervals around such roots with width
below $\tau_x$ to be good enough.  Clusters of zeros
of $\wit{f}$, as in Figure \ref{figc}, are a different kind
of beast, as the experience with other problems in
interval arithmetic has show \cite{Cluster1,ClusterN}.
By using the same tolerance $\tau_x$ for them we can
end up with literally thousands of small intervals
as candidates for containing roots, in the favorable
case in which the program actually ends and returns
something. We say this from our experience 
with our own  interval root finders as well as  
with interval root finders developed by other people. 

Our algorithm asks the user to also provide a ``cluster tolerance''
$\tau_c$, which could be something of the order of $\sqrt{\tau_x}$, and a tolerance
$\tau_w$ so that function values $w = \wfc{f}{t}$ with magnitude
below $\tau_w$ are deemed to be negligible. We take the 
liberty of increasing $\tau_w$ if we find it to be smaller than $16$ times
the maximum width of the intervals $\wfc{\wit{f}}{t}$ which compute along the execution
of algorithm. Contrary to what is proposed in \cite{HW}, we believe that
all tolerances should be absolute, and not relative (and users are free to
choose which approach suits their needs best.)
Using these tolerances, we can implement Algorithm \ref{alg:expand} which
``expands'' a zero $z \in \wit{x}$. By applying this algorithms to the point $z$ 
and the interval $\wit{x} = [\wu{x},\wo{x}]$ in Figure \ref{figc} 
we would identify the cluster $[\wu{c},\wo{c}]$ and submit the
intervals $[\wu{x},\wu{c}]$ and $[\wo{c},\wo{x}]$ to further
examination. Finally, we must say that the actual algorithms
for zero expansion in the \texttt{C++} code are more involved than
Algorithm \ref{alg:expand}, but the details are too cumbersome
to be presented here.

\begin{algorithm}[!h]
\caption{Zero expansion}
\label{alg:expand}
\begin{algorithmic}
\Procedure{expand\_zero}{$z$, $\wit{x}$}  

  \State $\wo{c} \gets z$
  \While{ $\wo{c} + \tau_c \leq \wo{x}$ }
		\State $\wit{w} \gets \wfc{\wit{f}}{\wo{c} + \tau_c}$ 
  	\State $\tau_w \gets \max \wset{16 \ \wfc{\wrm{wid}}{\wit{w}}, \tau_w}$
  	\If{ $\wu{w} \geq \tau_w$ or $\wo{w} \leq -\tau_w$ }
  		\State break
  	\EndIf
		\State $\wo{c} \gets \wo{c} + \tau_c$  
 \EndWhile
\State
  \State $\wu{c} \gets z$
  \While{ $\wu{c} - \tau_c \geq \wu{x}$ }
		\State $\wit{w} \gets \wfc{\wit{f}}{\wu{c} - \tau_c}$ 
  	\State $\tau_w \gets \max \wset{16  \ \wfc{\wrm{wid}}{\wit{w}}, \tau_w}$
  	\If{ $\wu{w} \geq \tau_w$ or $\wo{w} \leq -\tau_w$ }
  		\State break
  	\EndIf
		\State $\wu{c} \gets \wu{c} - \tau_c$  
 \EndWhile
\State
  \If{$\wu{c} > \wu{x}$ }
  \If{$\wo{c} < \wo{x}$ }
 		\State \Return $[\wu{x},\wu{c}], \ \ [\wu{c},\wo{c}],  \ \ [\wo{c},\wo{x}]$
  \Else
     \State \Return $[\wu{x},\wu{c}],\ \ [\wu{c},\wo{c}], \ \ \emptyset$
 \EndIf 
   \Else
  \If{$\wo{c} < \wo{x}$ }
  		\State \Return $\emptyset, \ \ [\wu{c},\wo{c}],  \ \ [\wo{c},\wo{x}]$
  \Else
			\State \Return $\emptyset, \ \ [\wu{c},\wo{c}],  \ \ \emptyset$
 \EndIf
  \EndIf

\EndProcedure
\end{algorithmic}
\end{algorithm}

\section{The modified Newton's method}
\label{sec_modify}
As mentioned in the introduction, the main tool for root finding is Newton's method
\begin{equation}
\label{newton2}
t_{k+1} = t_k - \wfc{f}{t_k}/\wdfc{f}{t_k}.
\end{equation}
This method was devised to find a ``point approximation'' to the root $r$, that
is, we for $\tilde{r}$ such that $\wabs{r - \tilde{r}}$ is small. 
It is our opinion that the goal of interval arithmetic
is a bit different: we look for a short interval $\wit{r}$ such that
$r \in \wit{r}$. Of course, in the end both goals amount to 
the same, but they suggest slightly different perspectives. 
We can rephrase the interval arithmetic goal as finding
a short interval $\wit{r} =[\wu{r},\wo{r}]$ such that 
$\wfc{f}{\wu{r}}$ and $\wfc{f}{\wo{r}}$ have opposite signs
(modulo degenerate cases in which $\wfc{f}{\wu{r}} = 0$ or $\wfc{f}{\wo{r}} = 0$.)
From this perspective, we believe that
we should modify the classic Newton step \wref{newton2} so that it produces
iterates $t_k$ such that the signs of $w_k := \wfc{f}{t_k}$ alternate,
in order to ensure that the interval
\[
\wit{r}_k := [ \wu{r}_k, \wo{r}_k] 
\whs{0.5} \wrm{for} \whs{0.5}
\wu{r}_k := \min \wset{t_k, t_{k+1}}
\whs{0.5} \wrm{and} \whs{0.5}
\wo{r}_k := \max \wset{t_k, t_{k+1}}
\] 
always contain a root. 

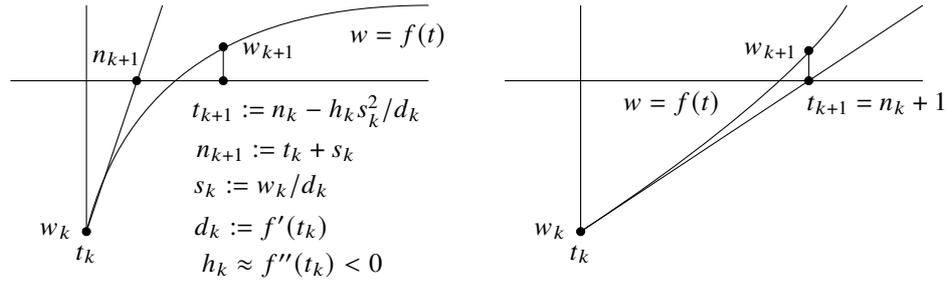
\begin{figure}[!h]
\begin{tikzpicture}[scale = 1]

\draw (0.0, 2.0)--( 5.5, 2.0);

\draw (1.0, 0.0)--( 1.0, 3.0);
\filldraw[black](1.0,0.0) circle(1.5pt);
\draw (1.0, 0.0) .. controls (1.5, 2.0) and (3.0, 3.0) .. (5.5,3.0);
\draw (1.0, 0.0)--( 2.0, 3.0);
\filldraw[black](1.66,2.0) circle(1.5pt);

\draw node at ( 0.6,  0.0){$w_k$};
\draw node at ( 1.0, -0.3){$t_k$};
\draw node at ( 1.4,  2.3){$n_{k+1}$};
\draw node at ( 5.1,  2.6){$w = \wfc{f}{t}$};
\filldraw[black](2.8, 2.0) circle(1.5pt);
\draw node at ( 3.9,  1.60){$t_{k + 1} := n_k - h_k s_k^2/d_k $};
\draw node at ( 3.5,  1.05){$n_{k+1} := t_k + s_k$};
\draw node at ( 3.3,  0.6){$s_k := w_k/d_k$};
\draw node at ( 3.3,  0.05){$d_k := \wdfc{f}{t_k}$};
\draw node at ( 3.7, -0.45){$h_k \approx \wdsfc{f}{t_k}< 0$};

\draw ( 2.8,  2.0)--(2.8, 2.4);
\draw node at ( 3.4,  2.4){$w_{k + 1}$};
\filldraw[black](2.8, 2.45) circle(1.5pt);

\draw (6.5, 2.0)--(12.0, 2.0);

\draw (7.5, 0.0)--( 7.5, 3.0);
\filldraw[black](7.5,0.0) circle(1.5pt);
\draw (7.5, 0.0) .. controls (10.5, 2.0) and (11.0, 3.0) .. (11.0,3.0);
\draw (7.5, 0.0)--(12.0, 3.0);
\filldraw[black](10.5, 2.0) circle(1.5pt);
\filldraw[black](10.5, 2.4) circle(1.5pt);
\draw node at ( 8.7,  1.7){$w = \wfc{f}{t}$};
\draw node at ( 7.1,  0.0){$w_k$};
\draw node at ( 7.5, -0.3){$t_k$};
\draw node at (11.4,  1.7){$t_{k + 1} = n_k+1$};
\draw node at (10.0,  2.4){$w_{k + 1}$};
\draw (10.5, 2.0)--(10.5, 2.4);

\end{tikzpicture}
\caption{The modified Newton step when $w_k = \wfc{f}{t_k} < 0$ and $d_k = \wdfc{f}{t_k} > 0$.
At left, $f$ is concave and we modify the Newton iterate $n_{k+1}$ in 
order to ensure that $w_{k+1}$ and $w_k$ have opposite signs. At right, $f$ is convex
and no modification is needed.}
\label{figi}
\end{figure}

A simple modification with this purpose is described in Figure \ref{figi}.
This Figure illustrates only the case in which
$w_k  < 0$ and $d_k > 0$, but cases in which
$w_k > 0$ or $d_k < 0$ are similar and can be reduced
to the case $w_k  < 0$ and $d_k > 0$ by replacing $f$
by $-f$ or $t$ by $-t$. A simple way to obtain a reasonable 
approximation $h_k$ for the second derivative 
mentioned in Figure \ref{figi} is to use the quotient
\[
h_k = \frac{d_{k-m} - d_k}{t_{k-m} - t_k},
\]
with data $d_{k-m}$ and $t_{k-m}$ from previous iterations.

Motivated by Figure \ref{figi}, we propose 
Algorithm \ref{alg:exact} below for finding a short
interval $\wit{r}$ containing the single root of $f$ in the interval
$\wit{x}$ using exact arithmetic, in the case in which
$\wfc{f}{\wu{x}} < 0 < \wfc{f}{\wo{x}}$ and $\wdfc{f}{t} > 0$ for 
$t \in \wit{x}$.
It is not difficult to prove that  Algorithm \ref{alg:exact}
has the same properties as the many other versions of Newton's
method that we find in the literature:
\begin{itemize}
\item When the model 
\begin{equation}
\label{quadratic}
\wfc{f}{t} \approx \wfc{f}{t_k} + \wdfc{f}{t_k} \wlr{t - t_k} + \wdsfc{f}{t} \wlr{t - t_k}^2/2
\end{equation}
is inaccurate for $t \in \wit{x}$, we fallback into a bisection step. This guarantees that 
eventually the interval $\wit{x}$ will become very short, and the analysis in the next item will
apply
\item When the quadratic model \wref{quadratic} is accurate in $\wit{x}$, 
and the interval $\wit{x}$ is very short,
$h$ will be a good approximation of $\wdsfc{f}{t_k}$ and our modification in the Newton step
will ensure that bisection will stop being necessary, the signs of the $w_k$ will alternate 
and the iterates will converge at the same quadratic rate as the classic method.
\end{itemize}

\begin{algorithm}[!h]
\caption{Find a short interval $\wit{r}$ containing the root $r$ of $f$ with $\wdfc{f}{t} > 0$.}
\label{alg:exact}
\begin{algorithmic}
\Procedure{exact\_newton}{$\wit{x}$, $f$, \ $\tau_x$, $\tau_w$}  
	\State  $\wu{w} \gets \wfc{f}{\wu{x}}, \ \ \wo{w} \gets \wfc{f}{\wo{x}}, \ \ \wu{d}\gets \wrm{nan}, \ \ \wo{d} \gets \wrm{nan}, \ \ d \gets \rm{nan}, \ \ t_d \gets \wrm{nan}, \ \ h \gets 0$

\If{ $\wu{w}> 0$ or $\wo{w} < 0$ }
	\State	\Return $\emptyset$
	\EndIf
		
  \State
	\State regular case:
  \If{ $\wlr{\wo{x}  - \wu{x} < \tau_x} \ $  or  $\ \wlr{\wo{w} - \wu{w} < \tau_w}$ }
		\State \Return $\wit{x}$  
  \EndIf	
	
  \If{ $-\wu{w} \leq \wo{w}$ }
  \If{ $!\wfc{\wrm{isnan}}{\wu{d}}$ }
  	\State goto bissection
  \EndIf
 	\State $t_k \gets \wu{x}, \ \ d_k \gets \wdfc{f}{t_k}$ 
 	\If{ $!\wfc{\wrm{isnan}}{d}$ and $t_d \neq t_k$}
 	   \State $h \gets \wlr{d - d_k}/\wlr{t_d - t_k}, \ \ d \gets d_k, \ \ t_d \gets t_k$
 	\EndIf
 	\State $\wu{d} \gets d_k, \ \ d \gets d_k, \ \ t_d \gets t_k, \ \ s_k \gets -\wu{w}/d_k$
  	\If{ $h < 0$ }
  		\State $s_k \gets s_k - h s_k^2/d_k$
  	\EndIf
  	\If{ $2 s_k > \wo{x} - \wu{x}$ }
  	  \State goto bissection
  	\EndIf
  	 \State $t_{k+1} \gets t_k + s_k, \ \ w_{k+1} \gets \wfc{f}{t_{k+1}}$
     \If{ $w_{k+1} \geq 0$ }
        \State $\wo{w} \gets w_{k+1}, \ \ \wo{x} \gets t_{k+1}, \ \ \wo{d} \gets \wrm{nan}$ 
        \State goto regular case
     \Else
     		\State $\wu{w} \gets w_{k+1}, \ \ \wu{x} \gets t_{k+1}, \ \ \wu{d} \gets \wrm{nan}$
       \State goto bissection
    	\EndIf
	\Else
	\State the case  $-\wu{w} > \wo{w}$ is analogous to $-\wu{w} \leq \wo{w}$ 
  \EndIf		
 
\State
\State bissection:
	\State $t_{k+1} \gets \wlr{\wu{x} + \wo{x}}/2, \ \  w_{k+1} \gets \wfc{f}{t_{k+1}} $	
	\If{ $w_{k+1} > 0$ }
	\State $\wo{w}  \gets w_{k+1}, \ \ \wo{x} \gets t_{k+1}, \ \  \wo{d} \gets \wrm{nan}$
	\Else
   \State $\wu{w}  \gets w_{k+1}, \ \ \wu{x} \gets t_{k+1}, \ \ \wu{d} \gets \wrm{nan}$
	\EndIf
\State goto regular case

\EndProcedure
\end{algorithmic}
\end{algorithm}

\section{The Monotonic Method}
\label{sec_mono}
This Section presents an interval Newton's method for functions
with $\wdfc{f}{t} > \kappa > 0$ for $t \in \wit{x}$ (in order to
used for functions with $\wdfc{f}{t} < -\kappa < 0$ we can
simply replace $f$ by $-f$. The method is a simple adaptation of the
``point method'' presented in Section \ref{sec_modify}, taking
into account the fundamental differences between exact (point) arithmetic
and interval arithmetic. Among such differences we have 
the following ones:
\begin{itemize}
\item In interval arithmetic, for a strictly increasing function, we may
have that $0 \in \wfc{\wit{f}}{t}$ for many $t$'s. In exact arithmetic
this can only happen for one value of $t$.
\item  In interval arithmetic we have a natural measure for the error
in the the values of $f$, given by the maximum width of the 
intervals $\wfc{\wit{f}}{t_k}$. This allows us to correct a poor choice
of the tolerance $\tau_w$ by the user and dispense with the cluster
tolerance $\tau_c$ mentioned in Section \ref{sec_expect}.
\end{itemize}

As we mentioned in the introduction, the increasing case  is much simpler than
the general one. In this case, we do not need to compute neither $\wfc{\wit{f}}{\wit{x}_k}$
nor $\wdfc{\wit{f}}{\wit{x}_k}$. We only need the interval $\wfc{\wit{f}}{t_k}$ and 
the floating point number $\wdfc{f}{t_k}$, and this entails a considerable reduction
in the cost of the steps. The analysis of degenerate cases is also much
simpler in the increasing case. These ideas are explored in 
Algorithm \ref{alg:inter_inc}, which is presented after the bibliography. 
For brevity, we omit the \texttt{expand\_zero}
label in this code, but it is similar to the zero expansion algorithms
presented in Section \ref{sec_expect}.

\begin{algorithm}[!h]
\caption{Find a short interval $\wit{r}$ containing the root $r$ of $f$ with 
$\wdfc{\wu{\wit{f}}}{\wit{x}} \geq \kappa > 0$.}
\label{alg:inter_inc}
\begin{algorithmic}
\Procedure{increasing\_interval\_newton}{$\wit{x}$, $f$, \ $\tau_x$, $\tau_w$}  
	\State  $\wit{w}_- \gets \wfc{\wit{f}}{\wu{x}}, \ \ \wit{w}_+ \gets \wfc{\wit{f}}{\wo{x}}, \ \
		 \wu{d}\gets \wrm{nan}, \ \ \wo{d} \gets \wrm{nan}, \ \ d \gets \rm{nan}, \ \ t_d \gets \wrm{nan}, \ \ h \gets 0$
		
	\If{ $\wu{w}_- > 0$ or $\wo{w}_+ < 0$ }
		\State	\Return $\emptyset$
	\EndIf
		
	\State $\tau_w \gets \max \wset{16 \ \wfc{\wrm{wid}}{ \wit{w}_-}, 16 \ \wfc{\wrm{wid}}{ \wit{w}_+}, \tau_w}$
	
	\If{ $\wo{w}_-	 \geq 0} $
		\State $t_z \gets \wu{x}$ and goto expand zero
	\EndIf
		
	\If{ $\wu{w}_+	 \leq 0} $
		\State $t_z \gets \wo{x}$ and goto expand zero
	\EndIf
	
	\State
	\State regular case:
  \If{ $\wlr{\wo{x}  - \wu{x} < \tau_x} \ $  or  $\ \wlr{\wo{w}_+ - \wu{w}_- < \tau_w}$ }
		\State \Return $\wit{x}$  
  \EndIf	
	
  \If{ $-\wu{w}_- \leq \wo{w}_+$ }
 
    \If{ $!\wfc{\wrm{isnan}}{\wu{d}}$ }
  	   \State goto bissection
    \EndIf
    
 	  \State $t_k \gets \wu{x}, \ \ d_k \gets \max \wset{\kappa, \wdfc{f}{t_k}}$ 
  
	 	\If{ $!\wfc{\wrm{isnan}}{d}$ and $t_d \neq t_k$}
	     \State $h \gets \wlr{d - d_k}/\wlr{t_d - t_k}, \ \ d \gets d_k, \ \ t_d \gets t_k$
 		\EndIf
 	
   	\State $\wu{d} \gets d_k, \ \ d \gets d_k, \ \ t_d \gets t_k, \ \ s_k \gets -\wu{w}_-/d_k$
   	
  	\If{ $h < 0$ }
  		\State $s_k \gets s_k - h s_k^2/d_k$
  	\EndIf
  	
  	\If{ $2 s_k > \wo{x} - \wu{x}$ } 
  	\State goto bissection
  	\EndIf
	
	 \State $t_{k+1} \gets t_k + s_k, \ \ {\wit{w}_{k+1} \gets \wfc{\wit{f}}{t_{k+1}}, \ \
           \tau_w \gets \max \wset{16 \ \wfc{\wrm{wid}}{ \wit{w}_{k+1}}, \tau_w} }$

     \If{ $\wu{w}_{k+1} > 0$ }
        \State $\wit{w}_+ \gets \wit{w}_{k+1}, \ \ \wo{x} \gets t_{k+1}, \ \ \wo{d} \gets \wrm{nan}$ and goto regular case
     \Else
     		\If{ $\wo{w}_{k+1} < 0$ }
       		\State $\wu{w}_- \gets \wit{w}_{k+1}, \ \ \wu{x} \gets t_{k+1}, \ \ \wu{d} \gets \wrm{nan}$
          \State goto bissection	
     		\Else
     		  \State $t_z \gets t_{k+1}$ and goto expand zero
        \EndIf
      \EndIf

	\Else
	\State the case  $-\wu{w}_- > \wo{w}_+$ is analogous to $-\wu{w}_- \leq \wo{w}_+$ 
  \EndIf		

\State
\State bissection:
	\State $t_{k+1} \gets \wlr{\wu{x} + \wo{x}}/2, \ \  \wit{w}_{k+1} \gets \wfc{\wit{f}}{t_{k+1}} $	
	\If{ $\wu{w}_{k+1} > 0$ }
	\State $\wit{w}_+  \gets \wit{w}_{k+1}, \ \ \wo{x} \gets t_{k+1}, \ \ \wo{d} \gets \wrm{nan}$
	\Else
		\If{ $\wo{w}_{k+1} > 0$ }
		   \State $\wit{w}_-  \gets \wit{w}_{k+1}, \ \ \wu{x} \gets t_{k+1}, \ \ \wu{d} \gets \wrm{nan}$
		\Else 
		   \State $t_z \gets t_{k+1}$ and goto expand zero
		\EndIf
	\EndIf
\State goto regular case

\EndProcedure
\end{algorithmic}
\end{algorithm}

\section{Testing}
\label{sec_test}
In practice, it is quite difficult to implement the algorithm presented
here due to the need of much attention to detail. It is quite likely
that the first attempts at such an implementation will fail, 
as ours did. Even worse, it is also likely that the failure will 
be undetected, specially if the details involved in the actual
coding of the algorithm are considered to be of little relevance. As a result,
in order to have some assurance of the reliability of 
our code in practice it is necessary
to have a good suit of tests. Readers should not
underestimate this point (and should not underestimate
our remark that they should not underestimate this point,
and, recursively...)

It is also important to realize that a good test suit 
may reveal problems not only in the code but also
in the theory upon which it is based. A good example
of this point is the choice of the stopping criterion
presented in chapter 9 of \cite{HW}. When using that
stopping criterion for the simple polynomial
\begin{equation}
\label{bughw}
\wfc{f}{x} = \wlr{x - 1} \wlr{x-2}\wlr{x- 3} \wlr{x-4} \wlr{x-5}
\end{equation}
in the interval $\wit{x} = [1,5]$, we would obtain that
the solution set $\ws{R}$ is the whole interval $\wit{x}$, regardless
of the tolerances provided by the user, and a good set of
tests could detect this problem. The authors of
\cite{HW} do mention that their stopping criterion may
be problematic in some rare cases, and that users
would detect this kind of problem afterwards.
Readers may think that the polynomial in Equation \wref{bughw} is one 
of such cases, but we beg to differ: we believe that 
a proper algorithm should  get the roots
of a function as simple as the one in Equation
\wref{bughw} with high accuracy without user intervention, 
and if it does not then it must be fixed.
Asking users to analyze the result may not be 
an option for instance when our algorithm is being used
thousands of times as part of a larger routine. In this 
scenario, users will not have the chance
to look at the result of each call of the algorithm,
and they will need either to rely on the algorithm to provide good answers,
or to write their own code to do what the algorithm should
have done for them.

The example in Equation \wref{bughw} can be though as a mild version of 
Wilkinson's polynomial, and we believe that we can obtain a good
test suit by adding multiple roots to polynomials similar
to it. In our opinion, a robust root finder must not
choke when trying to find the roots of polynomials of the form
\begin{equation}
\label{cases}
\wfc{p_{d,\wv{e}}}{x} := \pm \prod_{i = -m}^{m} \wlr{x -i}^{e_i}
\whs{0.4} \wrm{for} \whs{0.4} x \in \wit{x} = [m - \wu{\delta}, \ m + \wo{\delta}]
\end{equation}
where $m$ and the exponents $e_i$ are integers and
\[
\wu{\delta},\wo{\delta} \in \wset{0,1}, \whs{0.4}
\sum_{i=-m}^m e_i = d > 0, \whs{0.4}
\whs{0.4} m > 0 \whs{0.4}
 \wrm{and} \whs{0.4}
e_i \geq 0 \whs{0.2} \wrm{for} \whs{0.2} i = -m, \dots, m. 
\]
If $d$ and $m$ are not large then 
the coefficients of $p_{d,\wv{e}}$ can be computed
exactly in floating point arithmetic. For instance, 	when
$d = 20$ and $m = 5$ the coefficients of $p_{d,\wv{e}}$ can be 
computed exactly with the usual 
IEEEE 754 double precision arithmetic. There are 
\begin{equation}
\label{ndm}
n_{d,m} := 8 \ \binom{2 m + d}{d}
\end{equation}
elements in the family of polynomials in Equation \wref{cases}.
For $m = 5$ we have that 
\[
n_5 = \sum_{d = 1}^{20} n_{d,m} \approx 700 \ \wrm{million},
\]
and for tolerances  
$\tau_x = \tau_w = 10^{-6}$ and $\tau_c = 0.001$ with our code we can test all these $700$ million 
cases in a couple of days on our desktop machine, which has an AMD® Ryzen 7 2700x eight-core 
processor with 64MB of RAM (see the discussion about tolerances in Section 
\ref{sec_expect}.)

The  family of functions in Equation \wref{cases} has 
the following features, which make it a good test suit:
\begin{enumerate}
%%%%
\item By taking $\wu{\delta} = 0$ or $\wo{\delta} = 0$ we can check
how our code handles roots on the border of the interval $\wit{x}$.
This situation is a common source of bugs in root finders.
%%%%
\item Usually some $e_i$ will be greater than one, and the polynomials
$p_{d,\wv{e}}$ tend to have multiple roots. In fact, they may
have roots with multiplicity as high as $d$, or a couple of
roots with multiplicity $d/2$. These problems are quite challenging
and some root finders may return thousands of candidates intervals
for containing such multiple roots. 
%%%%
\item It is easy to write an algorithm to generate these polynomials 
by writing them in the factored form
\begin{equation}
\label{decomp}
\wfc{ p_{d,\wv{e}} }{x} = \pm q_1(-x) x^{e_0} q_2(x)
\end{equation}
where the $q_k$ are polynomials of the form
\[
\wfc{q_k}{x} = \prod_{i=1}^m \wlr{x - i}^{e_i}
\]
with degree at most $d$. For $m = 5$ and $d \leq 20$ there
are
\[
n_q = \sum_{j = 1}^{20} \binom{m + j - 1}{j} = 53129
\]
polynomials $q_k$ and we can keep them in a table
in memory, and let several threads build the polynomials
$p_{d,\wv{e}}$ from them using Equation \wref{decomp}.

\item When the coefficients $c_k$ of the expanded version of $p_{d,\wv{e}}$ 
\begin{equation}
\label{expand}
\wfc{p_{d,\wv{e}}}{x} = \sum_{k = 0}^d c_k x^k
\end{equation}
can be computed exactly, we know the exact roots of the polynomial
in Equation \wref{expand} and we can use this knowledge to evaluate
the robustness of our code with respect to the inaccuracies in the
evaluation of $p_{d,\wv{e}}$ and its derivatives using the expanded form.
In other words, we should test our code using Horner's method or
a Taylor model version of it to evaluate $p_{d,\wv{e}}$ 
using its expanded form \wref{expand}, because the evaluation of the
expanded form is usually much less accurate than the evaluation of the
product form in Equation \wref{cases}, and the purpose of testing
is to stress our code.
\end{enumerate}

Once we have chosen the family of functions and the intervals
on which will search for roots, we must decide which tolerances
to use for testing.  Our advice is to use both large and small
tolerances. Using large tolerances, like $0.1$, $0.01$ and $0.001$,
we can detect gross mistakes, like mistyping. 
We may miss some of these gross mistakes 
if we only use tiny tolerances, because with tiny tolerances
the part of the code affected by mistyping may
never be executed during the tests, or may only cause bugs which
are not severe enough to be detected
(we say this based on our personal experience
with our own blunders.) Tests with small tolerances are
necessary to access the algorithm's accuracy and to check
how much time it takes to find accurate roots.

\end{document}